\documentclass[12pt]{article}
\usepackage{amssymb}
\usepackage{latexsym}
\usepackage{amsmath,amscd}

\def\part#1{\frac{\partial\phantom{q}}{\partial#1}}

\newenvironment{rmk}{\begin{trivlist}\item[]{\bf Remark:} }
{\end{trivlist}}
\newenvironment{rmks}{\begin{trivlist}\item[]{\bf Remarks:} }
{\end{trivlist}}
\newenvironment{ex}{\begin{trivlist}\item[]{\bf Example:} }
{\end{trivlist}}
\newenvironment{prf}{\begin{trivlist}\item[]{\bf Proof:} }
{\hfill $\Box$ \end{trivlist}} 
\newtheorem{thm}{Theorem}

\newtheorem{prp}[thm]{Proposition}

\newcommand{\lie}[1]{\mathfrak{#1}}
\def\End{\mathop{\rm End}\nolimits}

\def\Hom{\mathop{\rm Hom}\nolimits}

\def\tr{\mathop{\rm tr}\nolimits}

\def\ad{\mathop{\rm ad}\nolimits}

\def\Diff{\mathop{\rm Diff}\nolimits}
\def\Ham{\mathop{\rm Ham}\nolimits}

\newcommand{\R}{\mathbf{R}}
\newcommand{\C}{\mathbf{C}}

\newcommand{\RP}{{\mathbf R}{\rm P}}
\newcommand{\PP}{{\rm P}}

\begin{document}
\title{A Teichm\"uller space for negatively curved surfaces}
\author{Nigel Hitchin\\{Mathematical Institute,
Woodstock Road,
Oxford, OX2 6GG}\\{hitchin@maths.ox.ac.uk}}

\maketitle
\begin{abstract}
We first describe the action of the fundamental group of a closed surface $\Sigma$ of variable negative curvature on the oriented geodesics in its universal covering in terms of a naturally-defined flat  connection whose holonomy lies in the group of Hamiltonian diffeomorphisms of $S^1\times \R$.  Consideration of the holonomy necessitates an extension from Riemannian to Finsler metrics. The second part of the paper follows the Higgs bundle approach to flat connections adapted to this infinite dimensional group and focuses on a family of  metrics, relying on a construction of O.Biquard, which is  parametrized by the 
infinite-dimensional space of CR functions on the unit  circle bundle of a hyperbolic surface. This generates an alternative approach to defining a  connection and offers the possibility of this vector space representing a moduli space which generalizes and includes the classical Teichm\"uller space. 
\end{abstract}

 \centerline{MSC classification:  53C22 53C60} 

\section{Introduction}
Classical Teichm\"uller space for a surface $\Sigma$ of genus $g>1$ has two manifestations: one is as a connected component of the space of homomorphisms  of the fundamental group $\pi_1(\Sigma)$ to $SL(2,\R)$, up to conjugation. The other  is  the space of metrics of constant negative curvature $-1$ up to diffeomorphism isotopic to the identity. The link is provided by uniformization -- the action of $\pi_1(\Sigma)$ on the universal covering of $\Sigma$, the upper half-plane $SL(2,\R)/SO(2)$. 

Now $SL(2,\R)$ acts also on the space of oriented geodesics in the hyperbolic plane. This is diffeomorphic to $S^1\times \R$ and the action preserves the standard symplectic form. Moreover this symplectic structure holds for the space of geodesics for the universal covering of {\it any} compact surface with  variable negative curvature and  gives representations of $\pi_1(\Sigma)$ in the infinite-dimensional group $\Ham(S^1\times \R)$ of Hamiltonian diffeomorphisms of $S^1\times \R$. These homomorphisms carry information about the metric --  by Crofton's formula they give the length of closed geodesics, and hence by Otal's theorem \cite{Otal} on the marked length spectrum, they define the isometry class. 

The first part of this paper gives a  construction of the symplectic connection which defines this representation and transforms naturally under  diffeomorphisms. In the second part we offer an approach to a moduli space of such connections which parallels the so-called Hitchin components of representations in $SL(n,\R)$ in terms of differentials of various degrees on $\Sigma$  with a fixed complex structure. 

From a differential geometric viewpoint  a  representation of $\pi_1(\Sigma)$  is provided by the holonomy of a flat connection $A$ on a bundle over $\Sigma$ with structure group $\Ham(S^1\times \R)$. This is   equivalent to a 4-manifold $\pi: M^4\rightarrow \Sigma$ fibring over $\Sigma$ with fibre $S^1\times \R$ and with a transverse foliation preserving the symplectic form along the fibres.  In Section \ref{sympr}  we construct a  flat connection by taking  $M=U\times \R$ where $U$ is the unit tangent bundle of a surface with Riemannian structure and negative curvature, and  introducing a natural exact decomposable  2-form. We identify the fibre over $x\in \Sigma$ with the geodesics in the universal covering by using the exponential map ${\mathrm {exp}}_x(tu)$ where $t\in \R$ and  $u\in U_x$, then  to $(u,t)\in S^1\times \R$ we associate the geodesic through ${\mathrm {exp}}_x(tu)$ intersecting  the geodesic $t\mapsto {\mathrm {exp}}_x(tu)$ with angle $\pi/2$. 

Having constructed a flat symplectic connection, we observe however that the holonomy  lies in a proper subgroup of $\Ham(S^1\times \R)$,  the diffeomorphisms which commute with the free antisymplectic involution $(\theta,t)\mapsto (\theta+\pi,-t)$ of $S^1\times \R$ whose quotient is the space of unoriented geodesics. This is the motivation for considering the more general case of {\it Finsler} surfaces of negative curvature in Section \ref{Fin}, where forward and backward geodesics are in general distinct.  After a brief introduction to Finsler surfaces, we adapt the Riemannian arguments to the new situation, which involves the consideration of different geodesic-like curves. 

The underlying motivation for this paper is to draw on the much-studied case of flat connections with holonomy in a noncompact Lie group, where generalizations of Teichm\"uller space were produced in \cite{Hit1}. In particular it was shown that for $SL(n,\R)$ there is a component of the moduli space which is diffeomorphic to a vector space and contains the classical Teichm\"uller space. From this viewpoint we choose to describe  
 $\Ham(S^1\times \R)$ as $SL(\infty,\R)$ (and the subgroup invariant by the involution as $Sp(\infty,\R)$). Given the notation  we might think of producing a moduli space which is a   ``classical limit" of the finite $n$ cases as $n$ tends to infinity, an idea which receives some support from  \cite{B2}. An alternative space with similar properties was introduced by Labourie  in \cite{Lab}. We put forward a concrete idea for this in the second part of the paper. 

Following a hint in the first order deformation theory in Section \ref{diffdef}, we introduce a family of Finsler metrics parametrized by an  infinite-dimensional vector space, each of which has an  associated 4-manifold describing a flat symplectic connection. 
 These metrics arise from a result of O.Biquard \cite{B1} following the ideas of \cite{Hit2} and we conjecture that these symplectic connections  are all equivalent to the connections defined earlier. 
 
  To be specific, from a sufficiently small CR function on the unit tangent bundle of a fixed hyperbolic surface Biquard's theorem produces a hyperk\"ahler 4-manifold defined on the interior of a disc bundle in $T^*\Sigma$ with a so-called fold singularity on the boundary.  There is an explicit one defined in terms of the hyperbolic metric and the full family consists  of deformations of this one. 
 The deformed boundary is then the  unit tangent bundle of a Finsler metric. Now the hyperk\"ahler metric has an analytic continuation to a neighbourhood of the exterior of the boundary as an indefinite metric of signature $(2,2)$ -- a {\it hypersymplectic} structure -- and this 4-manifold has a foliation by isotropic surfaces 
 which defines a flat symplectic connection whose leaves moreover intersect the unit tangent bundle in orbits  of the  geodesic flow. It is the global aspect of this feature which is missing to make the association with the first construction. 
  Nevertheless we calculate this hypersymplectic structure  in the case of quadratic CR functions  and recover the symplectic connection for any hyperbolic surface within this model. In the background is a discussion of  an involution on $SL(\infty,\R)$ like  the Cartan involution on $SL(n,\R)$ which provides an analogy with Corlette's approach to flat connections.

   The essential idea, for which this is a first step, is that a vector space -- the infinite-dimensional  space of CR functions -- can be identified with a component of the moduli space of  flat $SL(\infty,\R)$-connections, just as in the finite-dimensional case. 
  \subsection*{Acknowledgements}
The author wishes to thank  the EPSRC Programme Grant ``Symmetries and Correspondences"    and the Instituto de Ciencias Matem\'aticas, Madrid  for support during the preparation of this work. Thanks are also due to Olivier Biquard and Robert Bryant and to the anonymous referee of an earlier version of the paper. It was Jens Hoppe, many years ago, who first introduced me to the $n=\infty$ groups interpreted as symplectic diffeomorphisms \cite{JH}.  

\section{Symplectic connections}
\subsection{Covariant derivatives} \label{basic} 
A connection on a manifold $B$ with structure group the diffeomorphism group $\Diff(F)$ of another manifold $F$ is a fibration $\pi:M\rightarrow B$ with fibre diffeomorphic to $F$ and a horizontal distribution $H\subset TM$ which is transverse to the tangent bundle along the fibres and hence  $H\cong \pi^*TB$.  Writing $V$ for the tangent bundle along the fibres, the vertical bundle,  we have  a splitting $TM=V\oplus H$. The connection is flat if $H$ is integrable, defining a foliation. 

A connection $A$ on a principal bundle defines a covariant derivative on sections of the vector  bundle associated to the adjoint representation. In this case  given a  vector field $X$ on $B$ and a vertical vector field $Y$ the covariant derivative $\nabla_XY$ is the section of $V$  defined by lifting $X$ to a horizontal vector field $\tilde X$ and defining $\nabla_XY$ to be the vertical component of the Lie derivative $[\tilde X,Y]$. This extends in the usual way to a covariant exterior derivative 
$$d_A: C^{\infty}(M, \Lambda^pH^*\otimes V)\rightarrow C^{\infty}(M, \Lambda^{p+1}H^*\otimes V)$$
which defines a complex if the connection is flat. Since $H\cong \pi^*TB$ we could replace $H^*$ by $\pi^*T^*B$ and then the terms in the complex can be viewed as ``forms on $B$ with values in the Lie algebra of vector fields along the fibres".

Suppose now that $M$ is 4-dimensional,  $B=\Sigma$ is a surface and we are given a closed locally decomposable $2$-form $\alpha$ on $M$ which is nondegenerate on each fibre. Being decomposable (in four dimensions this is the condition $\alpha\wedge\alpha=0$) its annihilator defines a rank $2$ distribution $H$ and nondegeneracy on a fibre means that $H$ is transverse to the fibres. Since $\alpha$ is closed $H$ is integrable and this is a flat $\Diff(F)$-connection.

The $2$-form $\alpha$ gives each fibre a symplectic structure.  The Lie derivative  of $\alpha$ by a horizontal  vector field $\tilde X$ satisfies ${\mathcal L}_{\tilde X}\alpha=i_{\tilde X}d\alpha+d(i_{\tilde X}\alpha)=0$ because $i_{\tilde X}\alpha=0$ and $d\alpha=0$. It follows that   parallel translation along a curve maps one fibre to another symplectically. 

We are interested in the case where the fibre $F$ is $S^1\times \R$ which is not simply connected and so symplectic vector fields are not necessarily Hamiltonian.  If $\alpha$ is exact, then the connection can be lifted to a Hamiltonian one.  There is a choice  -- two Hamiltonians differ by a constant and so two Hamiltonian connections differ by the pullback of a one-form from $\Sigma$, and if the connections are flat then it is a closed form. 

For the symplectic connection, we can restrict the complex for the adjoint representation to Hamiltonian vector fields, in fact to Hamiltonian functions, giving
\begin{equation}
C^{\infty}(M)\stackrel {d_A}\rightarrow C^{\infty}(M,H^*)\stackrel {d_A}\rightarrow C^{\infty}(M,\Lambda^2H^*).
\label{complex}
\end{equation}
Here, if $f\in C^{\infty}(M)$, $d_Af$ is the component of $df$ in $H^*$ in the decomposition $T^*M=H^*\oplus V^*$.  

There are other descriptions:
\begin{itemize}
\item
 if $X$ is a vector field on $\Sigma$ and $\tilde X$ its horizontal lift, the covariant derivative of $f$ in the direction $X$ is $\tilde Xf$
 \item
$\alpha$ is  a non-vanishing section of 
$\Lambda^2V^*$ and embeds $H^*$ in $\Lambda^3T^*M$ then we can define $\alpha \wedge d_Af =d(f\alpha)$.
\item
Similarly if $a$  is a section of $H^*$, $\alpha\wedge d_Aa=d(\alpha \wedge a)$. From this it is clear that $d_A^2f=0$ since $d(\alpha\wedge df)=0$.
\end{itemize}

\subsection{Deformations} \label{def1}
Suppose we have a one-parameter family of exact locally decomposable $2$-forms $\alpha(t)=d\beta(t)$ then differentiating $\alpha\wedge\alpha=0$ at $t=0$ gives 
$$0=\alpha\wedge\dot\alpha=\alpha\wedge d\dot \beta=d(\alpha\wedge\dot\beta)$$
so the $H^*$-component $a$ of $\dot\beta$ satisfies $d_Aa=0$. 
\begin{prp} \label{inf} The section $a$ of $H^*\cong\pi^*T^*\Sigma$ represents the first order deformation of the symplectic connection defined by $\alpha$.
\end{prp} 
\begin{prf} Write $\dot\beta=a+b$ using the decomposition $T^*M=H^*\oplus V^*$, then $\dot\alpha=d(a+b)=(d_Ha, d_Va+d_Hb, d_Vb)$ is a section of $\Lambda^2H^*\oplus (H^*\otimes V^*)\oplus \Lambda^2 V^*$ and $d_Ha=d_Aa=0$ since $\alpha\wedge\dot\alpha=0$.

Using $\alpha$ as a nondegenerate 2-form on $V$ we have $$H^*\otimes V^*\cong H^*\otimes V\cong \Hom(H,V)=\Hom(H,TM/H)$$ and this component of $\dot\alpha$ is the first-order deformation of the horizontal distribution, namely the deformation of the $\Diff(F)$-connection. As above, a section of $H^*\otimes V$ is a one-form on $\Sigma$ with values in the Lie algebra of vector fields along the fibre, which is what a deformation of a connection is. 

The $\Lambda^2 V^*$ component $d_Vb$ is the variation of the symplectic form on the fibres defined by the restriction of $\dot\alpha$. We are supposed to be considering connections on a  fixed fibration by symplectic manifolds so, in the language of connections,  we need an infinitesimal gauge transformation to transform this to the standard one and apply this to the deformation to  get a variation relative to the same symplectic form. Now $b$ is a section of $V^*$ and so defines a 1-form on each fibre giving a  vector field $Y$ along the fibres with $i_Y\alpha=b$ and hence $d_Vb={\mathcal L}_Y\alpha$. This then is the variation of the symplectic structure -- transformed by $Y$. So $Y$ is the required infinitesimal gauge transformation and we must subtract $\nabla Y$, a section of $H^*\otimes V$,  from the variation of the connection to obtain a variation  for the {\it fixed} symplectic fibration. Identifying  $H^*\otimes V$ with  $H^*\otimes V^*$ using the symplectic structure, the variation of the connection is $d_Va+d_Hb$ and $\nabla Y$ identifies with $d_Hb$, so subtraction gives $d_Va$ as the deformation  of the symplectic connection. This is 
 clearly a Hamiltonian deformation defined by $a$, a section of $H^*$ in the complex.
\end{prf} 

\section{Riemannian surfaces}\label{sympr} 
\subsection{The symplectic connection} \label{connect} 
Let $\Sigma$ now be given a Riemannian metric $g$ of negative curvature. We shall describe a symplectic connection via a 4-manifold and a decomposable exact 2-form $\alpha$. First we approach the study of geodesics in a way which will extend easily to the Finsler case. 

Let $\pi: U\rightarrow \Sigma$ denote the unit circle bundle in $T\Sigma$. There are three canonical vector fields $X,Y,Z$ defined on $U$. The circle action on the fibres is given by $Z$, and the vector field giving the geodesic flow by $X$. The Levi-Civita connection defines a horizontal distribution and $X$ is tangent to this. The metric gives a complex structure $I$ on $\Sigma$ which induces a complex structure on the vector bundle $H$ and the third vector field is $Y=IX$ which is also horizontal. 

With the frame $X,Y,Z$ we have the dual basis $\omega,\theta, \eta$ and the following relations:

$\hskip 3cm[Z,X] =Y   \hskip 2cm d\omega =\eta\wedge\theta$
 
 $\hskip 3cm[Z,Y] = -X \hskip 1.62cm d\theta = -\eta\wedge\omega$
 
$\hskip 3cm [X,Y]= KZ \hskip 1.55cmd\eta= -K\omega\wedge \theta$

where $K$ is the Gaussian curvature. 

\begin{ex}
For the hyperbolic metric on the upper half-plane $y>0$ we have 
\begin{eqnarray*}
X&=&y\cos\phi\frac{\partial}{\partial x}+y\sin \phi\frac{\partial}{\partial y}-\cos\phi\frac{\partial}{\partial \phi}\\
 Y&=&-y\sin\phi\frac{\partial}{\partial x}+y\cos \phi\frac{\partial}{\partial y}+\sin\phi\frac{\partial}{\partial \phi}\\
  Z&=&\frac{\partial}{\partial \phi} 
\end{eqnarray*}
and 
\begin{eqnarray*}
\omega&=&\frac{1}{y}(\cos\phi \,dx+\sin\phi \,dy)\\
\theta&=&\frac{1}{y}(-\sin\phi \,dx+\cos\phi \,dy)\\
\eta&=& d\phi+\frac{dx}{y}.
\end{eqnarray*} 
\end{ex}

\begin{thm}\label{main}
Let $\psi_t$ be the flow generated by the vector field $Y$ on $U$ and define  $F(u,t)= \psi_t(u)$, then  
\begin{enumerate}
\item
the exact two-form $d(F^*\omega)$ on $M=U\times \R$ is decomposable and nondegenerate on the fibres of $\pi:U\times \R\rightarrow \Sigma$,
\item
the symplectic form on each fibre $S^1\times \R$ is diffeomorphic to the canonical symplectic form on $T^*S^1$.
\end{enumerate} 
\end{thm}
\begin{prf} 

Define one-forms on $U\times \R$ by  $\bar\omega=\psi^*_t\omega$ and similarly for $\theta,\eta$. Then the pulled-back form $F^*\omega$ of the theorem is $\bar \omega$.
We have 
$$d\bar \omega=\psi^*_t(\eta\wedge\theta+dt\wedge {\mathcal L}_Y\omega)=\bar\eta\wedge\bar\theta-dt\wedge \bar\eta=\bar\eta\wedge(\bar\theta+dt).$$
which is exact and decomposable. 

Evaluating these forms on $X$ and differentiating we obtain 
\begin{eqnarray*}
{\partial }\bar\eta(X)/{\partial t}&=&i_X\psi_t^*({\mathcal L}_Y\eta)=\bar K\bar\omega(X)\\
{\partial }\bar\theta(X)/{\partial t}&=&i_X\psi_t^*({\mathcal L}_Y\theta)=0\\
{\partial }\bar\omega(X)/{\partial t}&=&i_X\psi_t^*({\mathcal L}_Y\omega)=-\bar\eta(X)
\end{eqnarray*}
and similar results for $Y$,$Z$. Here the bar denotes also pulling back functions by $\psi_t$. 

We see in particular that 
  $\bar\omega$ evaluated on $X,Y,Z$ is a function $f$ which satisfies the equation $f''+\bar K f=0$. We can use these relations to express $\bar\omega,\bar\theta$ and $\bar\eta$  in terms of $\omega,\theta,\eta, dt$.
 
 At $t=0$, $\bar\omega(X)=\omega(X)=1$ and $\bar\omega(X)'=-\eta(X)=0$ so $\bar\omega(X)=f_1$ is the solution with initial condition   $f=1,f'=0$.  Similarly $ \bar\omega(Z)=f_2$  has initial conditions $f=0, f'=-1$.
For $\bar\omega(Y)$ we have   $f=0, f'=0$ so $\bar\omega(Y)$ vanishes identically and therefore 
$$\bar\omega=f_1\omega+f_2\eta$$

From the equation  $\bar\theta$ evaluated on $X,Y$ or $Z$ is constant and  by evaluation at $t=0$ we get  $\bar\theta=\theta.$
 
 As for $\bar\eta$, $\bar\eta(X)=-\bar\omega(X)'=-f_1'$, $\bar\eta(Y)'=0$ since $\bar\omega(Y)=0$ and the initial condition gives $\bar\eta(Y)=0$. Finally $\bar\eta(Z)=-\bar\omega(Z)'=-f_2'$ and 
 $$\bar\eta=-(f_1'\omega+f_2'\eta)$$

    Then 
\begin{equation}
d\bar\omega=\bar\eta\wedge(\bar\theta+dt)=-(f'_1\omega+f'_2\eta)\wedge (\theta+dt)
\label{alpha}
\end{equation}

Since $Z$ is tangential to the fibre $U_x\times \R$ over $x\in \Sigma$ then  $\omega$ and $\theta$ restrict to zero, so the restriction of  $F^*d\omega$ is  $-f'_2\eta\wedge dt$.

Now $f_2=\bar\omega(Z)$ satisfies $f''+\bar K f=0$ with initial conditions at $t=0$, $f=0, f'=-1$. Hence  for $t>0$ $f$ is negative up to and including the first point at which  $f'=0$. But at such a point $f''\ge 0$ and $f''+\bar K f=0$, which is impossible if $K$ and hence $\bar K$ is negative. A similar argument works for $t<0$. Then $f_2'$ is nonvanishing and the 2-form is nondegenerate restricted to a fibre.

The surface  $\Sigma$ is compact and so $\bar K$ is bounded away from zero which means that $f_2$ has exponential growth and is a diffeomorphism from $\R$ to $\R$. So we can reparametrize and write  the symplectic form as $df_2\wedge \eta$. Parametrizing the circle by $e^{i\phi}$ we have $Z=\partial/\partial \phi$ and $\eta =d\phi$ so the symplectic form is in standard form $ds\wedge d\phi$.
\end{prf}

\begin{rmk} The manifold $U\times \R$ in this setting can be regarded as the holonomy groupoid of the foliation of $U$ whose leaves are the orbits of $Y$. There is no good Hausdorff quotient of the foliation and the holonomy groupoid is in some sense the minimal desingularisation of the quotient. 
\end{rmk}

\begin{ex} For the hyperbolic plane $K=-1$ and $f_1,f_2$ are solutions of $f''-f=0$. With the given boundary conditions $f_1= \cosh(t), f_2=-\sinh(t)$ so that 
$$d\bar\omega=d(\cosh(t) \omega-\sinh(t) \eta)=(\theta+dt)\wedge (\sinh(t) \omega -\cosh(t) \eta).$$

We can compare this with the connection form for the known flat  $SL(2,\R)$ connection by realizing $S^1\times \R$ as a coadjoint orbit in the dual of the Lie algebra. Using the Killing form this is the quadric $u^2+v^2-w^2=1$ in $\R^3$. We take the symplectic form $dv\wedge du/w$ and  setting $w=\sinh t, u=\cosh t\cos\phi, v=\cosh t\sin \phi$ this   is $-\cosh t dt\wedge d\phi$, the form given above.  However, $u,v,w$ are also the Hamiltonian functions for generators $X,Y,Z$ of the Lie algebra of  $SL(2,\R)$  and  
using $\omega,\eta$ for the hyperbolic plane we  can then write the  form 
$\cosh(t) \omega-\sinh(t) \eta$ as $$\frac{1}{y}(\cosh (t)\cos\phi \,dx+\cosh (t)\sin\phi \,dy)-\sinh(t)(d\phi+\frac{dx}{y})=\frac{1}{y}(Xdx+Ydy)-Z(d\phi+\frac{dx}{y})$$
as a 1-form with values in ${\lie{sl}}(2,\R)$.
\end{ex} 
\subsection{The space of geodesics} 
Given a manifold of negative curvature, the exponential map identifies the universal covering $\tilde \Sigma$ with the tangent space at a point $x$. We can naturally parametrize geodesics in $\tilde \Sigma$ by taking $(u,t) \in U_x\times \R$,  and assigning to this data the geodesic which intersects the curve $s\mapsto \exp_x(su)$ orthogonally at $s=t$, taking the orientation $I\partial/\partial s$. If $t\ne 0$ then  $(u,t)$ is uniquely determined by the geodesic for any two such curves through $x$ will generate a geodesic triangle with angle sum greater than $\pi$, contradicting Gauss-Bonnet. At $t=0$ the direction $u$ determines  the oriented geodesic. This procedure gives an identification with $U_x\times \R=S^1\times \R$. The space of unoriented geodesics is diffeomorphic to the M\"obius band. 

The space of oriented geodesics $S$ has a  canonical  symplectic structure. A metric identifies the tangent bundle and the cotangent bundle and so the unit circle bundle $U$ sits naturally in both. In terms of our basis, the vector field $X$, being horizontal, is a section of $\pi^*T\Sigma$ and may be considered as the embedding of $U$ in $T\Sigma$. Similarly $\omega$, a section of $\pi^*T^*\Sigma$, is the embedding in the cotangent bundle. Equivalently 
 $\omega$ is the canonical one-form on $T^*\Sigma$ and $d\omega=\eta\wedge\theta$ the canonical symplectic form restricted to $U$.  Since $i_X\eta=0=i_X\theta$ the orbits of $X$ form the degeneracy foliation of $d\omega$ restricted to $U$. 
 
 If a closed 2-form has constant rank then it induces a symplectic form on the space of leaves of the foliation. In general this is a local statement but here, where  $U$ is the unit circle bundle of $\tilde \Sigma$, we have a global quotient diffeomorphic to $S^1\times \R$ which is therefore a symplectic surface.
 This surface is the space of orbits of the vector field $X$. The circular orbits of $Z$ project to circles which represent the points of $\tilde\Sigma$. 
  
 \begin{prp} \label{cover} The flat symplectic connection of Theorem \ref{main} on the universal covering $\tilde\Sigma$ is equivalent to the product $\tilde M=\tilde\Sigma\times S$ where $S$ is the space of oriented geodesics on $\tilde \Sigma$ with its canonical symplectic form $\alpha$.
 \end{prp} 
\begin{prf} The vector field $X$ on $U$ gave the geodesic flow but an integral curve of $Y=IX$ also projects to a geodesic on $\tilde \Sigma$ -- for each unit tangent vector $u$ it is the geodesic tangent to $Iu$. So the parametrization can be described in terms of the flow $\psi_t$ instead of the exponential map. 

Fix a point $x\in\tilde \Sigma$ and trivialize the $S^1\times \R$ bundle by parallel translation along radial geodesics $\gamma(s)$ through $x$ considered as  projections of integral curves  of $Y$. Since $d\bar\omega=\bar\eta\wedge(\bar\theta+dt)=-(f'_1\omega+f'_2\eta)\wedge (\theta+dt)$ the vector field $Y-\partial/\partial t$ on $M=U\times \R$ is horizontal, so parallel transport at time $s$  is the map
$$P(u,t)=(\psi_s(u),t-s)$$
from $(u,t)\in U_x\times \R$ to $U_{\gamma(s)}\times \R$. 

Each fibre $U_y\times \R$ is transverse to the foliation and the foliation intersects $U\times \{0\}$ in the geodesic flow given by the vector field $X$ since $d\bar\omega=d\omega$ at $t=0$. So the point $P(u,s)=(\psi_s(u),0)$ defines a geodesic through $\gamma(s)$. And since $Y=IX$ it is in the direction $-IY$. This establishes the identification of $U_x\times\R$ with the space of geodesics and parallel translation expresses $\tilde M\cong  \tilde \Sigma\times U_x\times \R$. 

 The closed form $d\bar\omega$  induces a symplectic form on the quotient by its degeneracy foliation and we now see that  this quotient   is well-defined and can  be identified with $U_x\times\R$. But this is also the quotient of $ U$ by the geodesic flow and there $d\bar\omega=d\omega$ so this is the canonical form on the space of geodesics. But  $d\bar\omega$ restricted to $U_y\times \R$ is the symplectic form preserved by the connection so  $\tilde M=\tilde\Sigma\times S$ where $S$ is the space of oriented geodesics on $\tilde \Sigma$. The identification depends of course on the basepoint $x$. 
\end{prf} 

\subsection{Holonomy}
Parallel translation along a curve from $x$ to $y$  provides a  diffeomorphism from $U_x\times \R$ to $U_y\times \R$ but it does not preserve this product structure, which is determined by the points $x,y$. We discuss here the relationship between points and geodesics under the action of the group $\Ham(S^1\times \R)$ of Hamiltonian diffeomorphisms.

   If we choose an identification of the space of geodesics with $S^1\times \R$, then 
 a point $x\in \tilde \Sigma$ defines a circle $C_x\subset S^1\times \R$, the geodesics passing through $x$. The circles $C_x,C_y$ corresponding to  points $x,y$ meet in two points, the geodesic joining $x$ and $y$ with its two orientations. They meet transversally since a Jacobi field (the function $f_2$) has at  most  one zero. 
The complement of $C_x\cup C_y\subset S^1\times \R$ has two bounded regions and the absolute value of the integral of the symplectic form over one of these gives twice the length of the geodesic segment $[x,y]$. This is a consequence of {\it Crofton's formula}. 

In our context we note that $\omega(X)=1$ and arc length along a geodesic is the parameter of the geodesic flow. Then the geodesic from $x$ to $y$ is the  image of  an  orbit of $X$ from $u\in U_x$ to $v\in U_y$ and with the opposite orientation from $-v\in U_y$ to $-u$ in $U_x$. Complete this to a 1-cycle by  semicircles in $U_x$, $U_y$ from $-u$ to $u$ and $v$ to $-v$ and integrate the 1-form $\omega$. Since $\omega(Z)=0$ this is twice the length but  the cycle bounds  a ``rectangular" 2-chain and by Stokes's theorem we obtain the integral of $d\omega$, the pull-back of the  symplectic form on $S^1\times \R$. The image of the rectangle is one of the bounded regions. The integral over the other region is the negative.

From Prop \ref{cover} the flat symplectic connection on the universal covering is a product $S^1\times \R\times \tilde\Sigma$ and since the orbit of $X$ through a point in the unit circle $U_x$ at $x\in \tilde \Sigma$  is horizontal, parallel translation preserves the point but transforms the circle through it.  A closed geodesic in $\Sigma$ lifts to a geodesic segment  in $\tilde \Sigma$ and the holonomy around that closed path defines a Hamiltonian diffeomorphism $g$.  Then the length of the geodesic is  half the symplectic area of one of the bounded regions in $C_x\cup g(C_x)$. 

Given any $g$ such that $g(S^1\times \{0\})$ intersects $S^1\times \{0\}$ transversally at two points, the symplectic area as above 
is an invariant under conjugation by Hamiltonian diffeomorphisms. A closed geodesic defines an element in $\pi_1(\Sigma)$ up to conjugation and if $g$ is the corresponding symplectomorphism  then  we see that the length of a closed geodesic  is defined by the conjugacy class of $g$, just as it is for hyperbolic geometry, with the finite-dimensional group $PSL(2,\R)$ acting symplectically. 
 
\subsection{ The groups $SL(\infty,\R)$ and  $Sp(\infty,\R)$}
The terminology $SU(\infty)$ for the symplectic diffeomorphisms of $S^2$ is  common in the physics literature (e.g.\cite{JH}). The reasoning is the decomposition of $C^{\infty}(S^2)/\R$, the Lie algebra of Hamiltonian diffeomorphisms, into irreducible representations of $SO(3)$. Each irreducible occurs with multiplicity one and this compares with the decomposition of $\lie{su}(n)$ under the action of the principal three-dimensional subgroup where all representations up to the $(2n-1)$-dimensional representation occur. As $n\rightarrow \infty$ the Poisson bracket also has an asymptotic form in terms of the finite-dimensional Lie bracket \cite{B2}. The integral of $f^2$ is then a substitute for the Killing form.

Instead of $S^2$, the hyperboloid $x^2+y^2-z^2=1$ is acted on by $SO(2,1)$ and this action is the same as $PSL(2,\R)$ on the space of geodesics on the hyperbolic plane, and for this reason we sometimes replace the notation $\Ham(S^1\times \R)$ by $SL(\infty, \R)$.  The parallel has its limits however, for although 
$PSL(2,\R)$ has a well-defined (indefinite) Killing form it does not have an extension to the Hamiltonian functions as in the case of $SU(\infty)$ -- and as we have seen the Hamiltonian functions for $PSL(2,\R)$ are not even square integrable. However, 
in the present context, a particular subgroup of $SL(\infty, \R)$ appears rather naturally. 

Although Riemannian metrics of negative curvature provide a large family of symplectic connections,  the holonomy has a special property.  Since a metric $g$ is a quadratic form on each  tangent space, $g(-u,-u)=g(u,u)$ and so $u\mapsto -u$ is an involution on the unit circle bundle $U$. This transforms the vector field $(X,Y,Z)$ to $(-X,-Y,Z)$. In particular the integral curve of  $Y$ for $t\ge 0$ is the integral curve of $-Y$ for $t\le 0$. It follows that the involution $\tau(u,t)=(-u,-t)$ on $U\times \R$ preserves the horizontal distribution of the symplectic connection. 

More analytically, the connection is defined by the form $d\bar \omega=d(f_1\omega+f_2\eta)$ and the relations $\tau^*\omega=-\omega, \tau^*\eta=\eta$ and  $\tau^*f_1=f_1, \tau^*f_2=-f_2$ give $\tau^*d\bar\omega=-d\bar\omega$. The involution acts in the fibres of $\pi:U\times \R\rightarrow \Sigma$ and since  the symplectic form  on the fibre is $f_2'dt\wedge\eta$ 
this action is antisymplectic.

 Define $Sp(\infty,\R)$ to  be the subgroup of the group $\Ham(S^1\times \R)$ consisting of elements which commute with the free antisymplectic involution $\sigma(\theta,t)=(\theta+\pi, -t)$
 (the notation for this group is motivated by the involution on $SL(2m,\R)$ whose fixed point set is $Sp(2m,\R)$).
With this definition we see that the holonomy of the symplectic connection for a Riemannian metric lies in the proper subgroup $Sp(\infty,\R)$. It acts on the quotient space $S^1\times \R/\tau$, the M\"obius band, which parametrizes unoriented geodesics. 
As a consequence of this observation we need to extend the type of geometry on the surface to obtain a larger space of representations.

\section{Finsler surfaces}\label{Fin} 
\subsection{Basic geometry}\label{finbas} 
There is a vast literature on Finsler geometry, but for surfaces we use the approach of Cartan, as described by Bryant in \cite{Br1}, \cite{Br2}. As originally conceived a Finsler manifold measures lengths of curves by a smoothly varying norm in each tangent space, but this is, by positive homogeneity, entirely determined by the unit circle bundle $U\subset T\Sigma$. The geometry  is then of a hypersurface in the tangent bundle such that $\pi:T\Sigma\rightarrow\Sigma$ is a submersion with fibres convex curves enclosing the origin. The circle bundle $U$ is the double covering of the projective tangent bundle $\PP(T)$ which in two dimensions is canonically $\PP(T^*)$ whose double cover is the unit cotangent vectors using the dual norm, so this same 3-manifold can be considered embedded in $T\Sigma$ or in $T^*\Sigma$. The degeneracy foliation for the restriction of the canonical symplectic form to 
$U\subset T^*\Sigma$  is tangent to the geodesic flow.

The circle bundle $U$ for a Finsler metric is no longer invariant under $u\mapsto -u$ and so in principle  may give more general symplectic connections. The lack of symmetry, or {\it reversibility}, means that there is an essential difference between forward geodesics and backward geodesics on $\Sigma$ not just one of orientation.

\begin{ex} The geodesics on a Riemannian manifold $(M,g)$ can be considered as the paths of a free particle, but one can also consider {\it magnetic geodesics}, describing the motion of a charged particle in a magnetic field given by a closed 2-form $F$. This is traditionally treated as  the Hamiltonian flow with respect to a new symplectic form $\Omega+\pi^*F$ where $\Omega$ is the canonical symplectic form on $T^*M$. If $F=d\alpha$ then this is equivalent to the geodesics for the Finsler metric whose unit sphere bundle is $U\subset T^*M$ translated by $-\alpha$. This is clearly not invariant by the involution. 
\end{ex}

As described in \cite{Br1} and \cite{Br2}, on a Finsler surface there is, just as in the Riemannian case,  a canonical basis $X,Y,Z$ for the tangent bundle  of $U$ and dual basis $\omega,\theta,\eta$ whose elements now satisfy the set of equations 

 $\hskip 2cm[Z,X] =Y   \hskip 5cm d\omega =\eta\wedge\theta$
 
 $\hskip 2cm[Z,Y] = -X+SY+CZ \hskip 2.4cm d\theta = -\eta\wedge(\omega+S\theta)$
 
$\hskip 2cm [X,Y]= KZ \hskip 4.55cmd\eta= -(K\omega +C\eta)\wedge \theta$

 Here $S,C,K$ are functions on $U$ and  we still call $K$ the (Gaussian) curvature. We follow here the notation in \cite{Br1} although $S,C$ are frequently called $I,J$. 
 
 The vector field $Z$ is tangential to the fibre, so $X$ and $Y$ span the horizontal space for a $\Diff(S^1)$-connection. Its curvature is $K$, or better $(KZ)\omega\wedge\theta$, a 2-form on $\Sigma$ with values in vector fields along the fibres. We shall consider Finsler surfaces where $K<0$. 
 The vector field $X$ generates the geodesic flow. 
 
It is a useful convention to use local coordinates $x,y$ on $\Sigma$ and $p$ on $U$  and then the Finsler structure is given by a function $L(x,y,p)$.  
In this scenario a curve $(x(t),y(t))$ lifts canonically to the tangent bundle as $x'\partial/\partial x+ y'\partial/\partial y$ and $p=y'/x'$. Then the length of a curve for $t\in [a,b]$ is defined to be 
$$\int_a^b x'L(x,y,p)dt=\int_{x(a)}^{x(b)} L(x,y,p)dx.$$
 In terms of the usual local coordinates $(u_1,u_2)\mapsto u_1\partial/\partial x+ u_2\partial/\partial y$ on $T\Sigma$ the unit circle bundle  $U$ is defined by  $u_1L(x,y,p)=1$ with   $p=u_2/u_1$. 

  The dual norm of $\xi_1dx+\xi_2dy$ is given by the supremum of $\xi_1u_1+\xi_2u_2$  where $u_1L(p)=1$,  or $\xi_1/L+\xi_2p/L$ so the maximum  occurs when 
  $-{L_p}/{L^2}(\xi_1+\xi_2)+\xi_2/L=0$
  giving the value $\xi_1/L+\xi_2p/L=\xi_2/L_p$. So $U$ is embedded in $T^*\Sigma$ as $$(\xi_1,\xi_2)=(L-pL_p,L_p).$$ The canonical one-form on $T^*\Sigma$ therefore  restricts to $(L-pL_p)dx+L_pdy$. This is the first element of the canonical basis 
 \begin{eqnarray*}
 \omega &=& Ldx+L_p(dy-pdx)\\
 \theta &=& \sqrt{LL_{pp}}(dy-pdx)\\
 \eta&=& \frac{dL_p-L_ydx}{\sqrt{LL_{pp}}}+\frac{A(L)(dy-pdx)}{{\sqrt{(LL_{pp})^3}}}
  \end{eqnarray*}
where (as in \cite{Br1} Prop 5.) $A(L)$ is a universal polynomial in $L$ and its derivatives determined by $\eta\wedge d\theta=0$. 

Note that $\theta$ vanishes on the canonical lift of a curve in $\Sigma$ to $U$: $p=dy/dx$,  in particular the geodesic flow gives the canonical lift of  a geodesic.

\begin{ex}
For the hyperbolic metric $L=\sqrt{1+p^2}/y$ or $L=1/y\cos\phi$ with $p=\tan \phi$, 
\end{ex} 
\subsection{$Y$-curves}

The vector field $X$  is again the geodesic flow but $Y$ now has a different role. In a Riemannian surface it was horizontal and orthogonal to $X$, but in a Finsler metric we don't have quite the same notion of angle. (There is an analogue called the {\it Landsberg angle} -- given two points $u,v\in U_x$ the angle is defined by the integral of $\eta$ from $u$ to $v$ normalized by the integral over $U_x$).

It is convenient for the purposes of this paper to give the name  $Y$-{\it curves} to the projections of the orbits of $Y$ on $\Sigma$. A recent paper \cite{ISS} calls them $N$-{\it parallels}.

Orthogonality and the role of $Y$ can be seen by  considering the variational problem of the shortest distance from a fixed point to a curve $y=f(x)$ (we shall use local coordinates here for convenience). If $\dot g(x)$ is the variation of the curve $y=g(x)$, the initial point $(x,y)=(a,g(a))$ is fixed so $\dot g(a)=0$ but the upper end-point $b=z$ has variation $\dot z$ so since $g(z)=f(z)$ we have $\dot g(z)+g'(z)\dot z=f'(z)\dot z$. Then the variation of the integral is 
$$ L\dot z+\int_a^z L_y\dot g+L_p\dot g'dx= L\dot z+L_p\dot g(z)+\int_a^z (L_y-(L_p)')\dot g dx$$
integrating by parts. The integral term gives the usual Euler-Lagrange equations, so $y=g(x)$  is a geodesic, and the vanishing of the first term is 
$$0=L\dot z+L_p\dot g(z)= (L+L_p(f'-g'))\dot z.$$
 At the point $m=(x,y,p)=(x,g(x),g'(x))$  the horizontal  subspace of $T_mU$ spanned by $X$ and $Y$ is isomorphic to the tangent space of $\tilde\Sigma$ at $(x,g(x))$ and $X$ is the tangent vector of the geodesic $y=g(x)$. Restricting the form  $\omega= Ldx+L_p(dy-pdx)$  to the tangent space of the curve $y=f(x)$ gives $(L+L_p(f'-g'))$ which vanishes in the current situation. But since $\omega(X)=1,\omega(Y)=0$ this is the direction of $Y$.

It is in this sense that  $X$ and $Y$ are orthogonal: the shortest path to a curve meets it orthogonally. In local coordinates $q=p-L/L_p$ is the relation. 

The integral curves of $Y$  no longer project to geodesics but they do share many properties which will become relevant later. 
The first concerns their deformations.

 First we revisit the space of geodesics and its symplectic structure: a tangent vector to the space of  orbits of $X$ is a section of the normal bundle which is invariant by the action of $X$. Using the canonical basis we can  write this section as $aY+bZ$ and $X$-invariance gives $(Xa)Y+aKZ+(Xb)Z-bY=0$ modulo $X$. Hence $b=Xa, Xb=-aK$ and $a$ satisfies the equation $X^2a+Ka=0$, or $a''+Ka=0$ along the flow. This equation for Jacobi fields is a self-adjoint ODE and the Wronskian of two sections is the symplectic form on the space of geodesics.

If we do the same for orbits of $Y$, a section of the normal bundle is $aX+bZ$ and invariance gives $(Ya)X-aKZ+(Yb)Z+b(X-SY-CZ)=0$ modulo $Y$ so in this case we obtain $a''-Ca'+Ka=0$ which is not self-adjoint unless $C=0$. 

For the sake of completeness we do the same for $Z$, whose orbit space is of course $\Sigma$. In this case $aX+bY$ is normal and invariance yields 
$(Za)X+aY+(Zb)Y+b(-X+SY+CZ)=0$ modulo $Z$ which is $Za=b, Zb+Sb=-a$ so $a''-Sa'+a=0$ for $a=a(\phi),b=b(\phi)$ with $\phi$ a parameter with  $Z=\partial/\partial \phi$. Each tangent vector at a point $x$ in $\Sigma$ defines a deformation of the fibre  and this  provides more information than is apparent from the equation $a''-Sa'+a=0$. The 2-dimensional tangent space to $\Sigma$  gives  two independent periodic solutions of this equation. Furthermore, since they define  lines through the origin in the tangent space at $x$, they have two simple zeros on $U_x$, like $\sin \phi,\cos\phi$.

This formalism is useful for understanding the $\Diff(S^1)$-connection associated to a Finsler surface. Given a curve $\gamma(t)$ 
in $\Sigma$ then the tangent vector field $\dot\gamma(t)$ along the curve can be lifted horizontally to a vector field $W$ in $U$ over $\gamma$. Integrating to a flow gives the holonomy from the point $x$ of $\gamma$ to another point $y$ as a diffeomorphism from $U_x$ to $U_y$. 

At a point $x\in \gamma$ we have two distinguished points in $U_x$: where $a=1,b=0$, the unit tangent vector at $x$ and where $a=0,b=1$. The latter we call the unit normal, from the discussion of the direction $Y$ above.
A geodesic lifts horizontally to an orbit of $X$, so that $b=0,a=1$ along the curve, meaning the unit tangent vector is preserved by parallel translation. 
Similarly a  $Y$-curve lifts horizontally to an orbit of $Y$ where $a=1,b=0$ along the curve, so parallel translation preserves the unit normal -- this is the motivation for the name  ``$N$-parallel".

\begin{rmk} As we observed, duality embeds $U$ in $T^*\Sigma$ but in two dimensions $T^*\cong T\otimes \Lambda^2T^*$ and $\theta\wedge\omega$ trivializes $\pi^*\Lambda^2T^*$ on $U$, so with this identification we obtain a new Finsler metric $\tilde L$. The vector field $Y$ now satisfies $\tilde \omega(Y)=1, \tilde \theta(Y)=0$ but is not a geodesic flow because, as remarked above, its deformation space lacks a natural symplectic form. 
\end{rmk}

We need now to revisit the issues we dealt with for Riemannian surfaces in the Finsler context. As described in \cite{Fang} the geometry on the universal covering of a compact negatively curved Finsler manifold has good properties -- the Cartan-Hadamard theorem holds, we have forward and backward geodesic completeness, convexity for open balls etc. There is unfortunately less discussion in the literature of properties of the $Y$-curves.
\subsection{The symplectic connection} 
The proof of Theorem 2 and Proposition 3 for negatively curved Finsler surfaces follows the same lines as the Riemannian case.  The first difference is that  the 1-forms  $\bar\omega$ etc. evaluated on $X,Y,Z$ now give  functions $f$ which satisfy the equation  $f''-\bar Cf'+\bar K f=0$. We introduce functions $f_1,f_2$ with the same initial conditions and 
then  $\bar \omega = f_1\omega+f_2\eta$ and  $\bar\eta  
  =-(f_1'\omega+f_2'\eta).$  
   However, now   $\bar\theta=c\omega+\theta+e\eta$  for functions $c,e$ with $\bar\theta'=-\bar S\bar\eta.$ The closed 2-form giving the connection becomes 
   $$-(f_1'\omega+f_2'\eta)\wedge (\theta+dt+c\omega+e\eta).$$
   Then $d\bar\omega=-f_2'dt\wedge \eta$ on each fibre $U_x\times \R$ and we need to prove that this is nondegenerate. The new equation is 
   $f''-\bar Cf'+\bar K f=0$. At $t=0$, $f_2=0, f_2'=-1$ and so for $t>0$ $f_2$ is negative up to and including the first point at which  $f_2'=0$. But at such a point $f_2'' =-\bar K f_2 \ge 0$ which is impossible if $K$ and hence $\bar K$ is negative. A similar argument works for $t<0$ and so $f_2'\ne 0$ for all $t$. 
   
  To carry this further,  $f_2'<0$ so $f_2$ is negative and decreasing. If it is bounded below as $t\rightarrow \infty$, say $f_2\ge -m$,  then $f_2'\rightarrow 0$. But  since $\bar K\le - M$, the equation  $f_2''=\bar C f_2'-\bar K f_2$ implies $f_2''\ge  -Mm$ which contradicts $f_2'\rightarrow 0$.
  It follows that $f_2$ is unbounded and $-f_2'dt\wedge \eta$ is diffeomorphic to the standard symplectic form on $S^1\times \R$ and we have a symplectic connection with holonomy in $SL(\infty,\R)$.

   To relate the symplectic connection to the action of $\pi_1(\Sigma)$ on geodesics we note first that  there is a map, as in the Riemannian case, from each fibre $U_x\times \R$ to the space of geodesics on the universal covering: at each point on a $Y$-curve through $x$ we take the orthogonal (forward) geodesic. We need to show that every geodesic is found this way and that the $Y$-curve through $x$ to it is unique.

    Given a geodesic $\gamma$, let $\gamma(s)$ be the canonical lift to $U$ and define $F(s,t)\in \tilde \Sigma$  to be the projection of the flow $\psi_t$ of $Y$ through $\gamma(s)$.

    Consider the derivative of $F:\R^2\rightarrow \tilde\Sigma$. As $s$ varies along $\gamma$ we have a 1-parameter family of $Y$-curves and hence at a point $\gamma(s)$, the tangent vector $\partial/\partial s$ gives a section of the normal bundle in $U$ of the $Y$-curve through $\gamma(s)$. As above this is of the form $aX+bZ$ where $Ya+b=0$ and $Y^2a-CYa+Ka=0$, or in terms of the parameter $t$ along the $Y$-curve $a'=-b$ and $ a''-Ca'+Ka=0$. The $Y$-curves meet the curve $\gamma(s)$ orthogonally at $t=0$ which means that $b=0$ there and since $X=\partial/\partial s$  we have $a=f_1(t)$, the solution with initial condition $f_1(0)=1, f_1'(0)=0$. By a similar argument to the discussion of $f_2$, $f_1\ge 1$ and is unbounded above. 
    
    Using the derivative $D\pi$ of the projection $\pi:U\rightarrow \tilde \Sigma$, we use $(D\pi(X), D\pi(Y))$ as a basis for $T\Sigma$  and the derivative of $F$ is 
    $$DF(\dot s, \dot t)=\dot s f_1(t,s)D\pi(X)+\dot t D\pi(Y).$$
    Since $f_1$ is  bounded below by $1$  this is invertible and hence $F$ is a local diffeomorphism from $\R^2$ to $\tilde\Sigma$. 
    
    Now $X$, the geodesic flow, has norm 1 and the norm of $Y$ is uniformly bounded above  and below by compactness of $\Sigma$, thus  the inverse of $DF$, namely $(X,Y)\mapsto (f_1^{-1},1)$, has  bounded norm. From Hadamard's global inverse theorem $F$  is a global diffeomorphism from $\R^2$ to $\tilde\Sigma$. It follows first of all that $F$ is surjective and so through any point $x\in \tilde\Sigma$ there exists a $Y$-curve orthogonal to the given geodesic, and secondly that curve is unique. 
    
      \begin{rmk} For a Riemannian metric the length of $X$ and $Y$ are $1$ for whatever the point in a fibre of $U$ so the metric pulled back by $F(s,t)$ is 
   $$f_1^2ds^2+dt^2= G(s,t) ds^2+dt^2.$$
  Elementary texts in differential geometry show that the Gaussian curvature is $K=-G^{-1/2}(G^{1/2})_{tt}$. In our case the equation satisfied by $f_1$ gives this since $K=-f_1''/f_1$.
     \end{rmk} 

\section{Deformations} \label{diffdef}
In this section we shall consider first order deformations of the flat symplectic connection on a Finsler surface, following   the general approach in Section \ref{def1}. 

A diffeomorphism of $\Sigma$ connected to the identity acts trivially on the fundamental group, hence the action of a vector field on $\Sigma$ should give us a first order variation of the symplectic connection of the form $d_Ah$  for some function $h$ on $M=U\times \R$. Regarding $U$ again as the double covering of the projective bundle $\PP(T)$, there is a natural lift of the action of a diffeomorphism on $\Sigma$  to $\PP(T)$, preserving the fibres. It is characterized by preserving both the fibration and the contact distribution. It follows that a vector field $W$ on $\Sigma$ induces a vector field $\tilde W$ on $U$. 

\begin{prp} The vector field $W$ on $\Sigma$ defines the variation $d_Ah$ of the symplectic connection where $h=i_{\tilde W}(f_1\omega+f_2\eta)$.
\end{prp} 
\begin{prf}
The connection is given by  $\alpha=d(f_1\omega+f_2\eta)=d\beta$ where $\beta$ is naturally defined by the Finsler metric so from Section \ref{basic}  the variation of the connection  is defined by the $H^*$ component of ${\mathcal L}_{\tilde W}\beta$. This in turn is uniquely determined by the 3-form $\alpha\wedge {\mathcal L}_{\tilde W}\beta$.
Consider $\alpha\wedge d(i_{\tilde W}\beta)=\alpha\wedge ({\mathcal L}_{\tilde W}\beta-i_{\tilde W}d\beta)$. Now $\alpha\wedge\alpha=0$ so taking the interior product with $\tilde W$ gives $0=\alpha\wedge i_{\tilde W}\alpha=\alpha\wedge i_{\tilde W}d\beta$ and  therefore $\alpha\wedge d(i_{\tilde W}\beta)=\alpha\wedge {\mathcal L}_{\tilde W}\beta$.

Since $\alpha\wedge dh=\alpha\wedge d_Ah$ for any function we see that $h=i_{\tilde W}(f_1\omega+f_2\eta)$ is the required function, showing that the cohomology class in the complex (\ref{complex}) defined by an infinitesimal diffeomorphism of $\Sigma$ is trivial, as expected.
\end{prf} 

There is a distinguished class of deformations for a {\it hyperbolic} surface, which we examine next.  Firstly, recall that the connection is defined by the 2-form 
$$\alpha=d(\cosh(t)\omega-\sinh(t) \eta)=(\theta+dt)\wedge(\sinh(t)\omega-\cosh(t)\eta).$$  
A general first order deformation  is then given by a section $a=a_1\omega+a_2\theta$ of $H^*$ with $d_Aa=0$ or equivalently $\alpha\wedge da=0$. This gives the equation
\begin{equation}
(Xa_2-Ya_1)+(Za_2+a_1)\tanh(t)+a_1'=0.
\label{DA}
\end{equation} 
Trivial deformation classes are of the form $d_Af=(Xf+\tanh(t) \,Zf)\omega +(Yf-f')\theta$.

Any Riemannian surface  $\Sigma$ has a complex structure on each tangent space   and this gives the unit tangent bundle $U$ a CR structure. Since $Y=IX$,  a complex function  $h$ is a CR function if  $(X+iY)h=0$.  Moreover the vector field $Z$ generates a circle action and 
 $[Z,X+iY]=-i(X+iY)$, so the space of CR functions is preserved.  
 
 \begin{prp} \label{CRprop} Let $h=h_1+ih_2$ be a CR function on the unit tangent bundle of a hyperbolic surface  with $Zh=imh$ and $m>0$. Then a non-zero solution to the deformation equation (\ref{DA}) is provided by 
 $$a_1-ia_2=h \cosh^{m-1}(t)$$
  \end{prp}
  \begin{prf} The formula gives $Z(a_1-ia_2)=im(a_1-ia_2)$ so $Za_2=-ma_1$. But $(X+iY)h=0$ implies $Xa_2-Ya_1=0$ so equation (\ref{DA}) holds if $(m-1)a_1\tanh(t)=a_1'$. 
  
  To see that this is nonzero if $m\ge 0$, observe that we have 
   the function $h$ transforming under the circle action as  $h\mapsto e^{im\phi}h$. Now $\omega$ embeds $U$ as the unit cotangent bundle with $\omega$  the restriction of the canonical one-form and $\omega+i\theta$ is the canonical holomorphic one-form on $T^*\Sigma$. With our conventions $\mathcal{L}_Z(\omega+i\theta)=-i(\omega +i\theta)$ as $\omega+i\theta=y^{-1}e^{-i\phi}(dx+idy)$.  
  
  We can also interpret $\omega+i\theta$ as the tautological section $s$ of $\pi^*K$, the pull-back of the canonical bundle to its total space, then $s^m$ is a section of $\pi^*K^m$ such that $\mathcal{L}_Zs^m=-ims^m$. Then $$\mathcal{L}_Zhs^m=(im-im)hs^m=0.$$ This means that $hs^m$ is invariant therefore  the pull-back of a section of $K^m$. Moreover  $(X+iY)h=0$ implies that this section is holomorphic. These are non-zero for $m\ge 0$ and form a vector space $H^0(\Sigma, K^m)$ of dimension $(2m-1)(g-1)$.
 \end{prf} 
 
\begin{rmk}
 When $m=1$  $a_1-ia_2$ is independent of $t$ and the deformation consists of  adding the real part of a holomorphic, and therefore closed, 1-form to $\cosh(t)\omega-\sinh(t)\eta$, a magnetic deformation. It does not change the actual connection defined by the 2-form $\alpha$. 
\end{rmk} 

We have identified here  distinguished directions to deform the symplectic connection of a hyperbolic surface, moreover the introduction of holomorphic sections of $K^m$ suggests an analogy with the component of the moduli space of $SL(n,\R)$ connections which generalizes Teichm\"uller space. As in  \cite{Hit1} this is represented by the direct sum $H^0(\Sigma, K^2)\oplus \cdots \oplus H^0(\Sigma, K^n)$. Moreover the analogous space for $Sp(2n,\R)$ consists of sections of the even powers. Could it be that the moduli space for $SL(\infty,\R)$ is isomorphic to the vector space of all CR functions with a suitable topology, and in particular be contractible? We certainly know from Ricci flow  \cite{Ham} that Teichm\"uller space is a deformation retract of the space of all negatively curved Riemannian surfaces. This issue we take up in the second part of the paper where we produce a symplectic connection in a rather different way. 

\section{$SL(\infty,\R)$ Higgs bundles}

\subsection{Comparison with $SL(n,\R)$}
We recall here how the theory of Higgs bundles produces a generalization of Teichm\"uller space for $SL(n,\R)$. 

Given an irreducible representation of $\pi_1(\Sigma)$ in $SL(n,\R)$ and a complex structure on $\Sigma$, the theorem of Corlette \cite{Cor} produces an equivariant harmonic map from the universal covering $\tilde \Sigma$ to the symmetric space $SL(n,\R)/SO(n)$. The reduction of structure group to $SO(n)$ gives a (non-flat)  connection $A$ over a principal $SO(n)$ bundle $P$ over $\Sigma$ and the flat $SL(n,\R)$ connection can be written as $\nabla_A+\phi$ where $\phi\in \Omega^1(\Sigma, \ad P)$ is symmetric. Using the complex structure on $\Sigma$, we can write $\phi=\Phi+\Phi^*$ where $\Phi\in H^0(\Sigma, \ad P^c\otimes K)$, the Higgs field, is a holomorphic section of the adjoint bundle tensored with the canonical bundle.  

Conversely, by work of Simpson (and the author for $n=2$), the pair of a holomorphic $SO(n,\C)$ bundle and Higgs field $\Phi$ satisfying a stability condition yields a canonical metric on the bundle such that the curvature $F_A$ of the associated $SO(n)$ connection satisfies the equation $F_A+[\Phi,\Phi^*]=0$ and $\nabla_A+\Phi+\Phi^*$ is then a flat $SL(n,\R)$ connection. This holomorphic point of view gives for each solution sections of $K^m$ for $2\le m\le n$ as the coefficients of the characteristic polynomial $\det(x-\Phi)$. 

When $n=2$  take the vector bundle $V=K^{-1/2}\oplus K^{1/2}$ and the Higgs field $\Phi(u,v)=(v,qu)$ where $q\in H^0(\Sigma,K^2)$. Then the determinant is $-q$ and this invariant
 identifies a component of the space of $SL(2,\R)$ representations with the vector space  $H^0(\Sigma,K^2)$. This is Teichm\"uller space but not with the natural complex structure. The distinguished origin $q=0$ is the uniformizing representation for the given complex structure. 
 
 For the higher Teichm\"uller spaces of representations into $SL(n,\R)$ we take symmetric powers of $K^{-1/2}\oplus K^{1/2}$:  for  $n=2m+1$ the vector bundle 
 $V=K^{-m}\oplus K^{1-m}\oplus \dots \oplus K^{m}$
 and for $n=2m$ we have 
 $V=K^{-(2m-1)/2}\oplus K^{1-(2m-1)/2}\oplus \dots \oplus K^{(2m-1)/2}.$
The pairing of $K^{\pm \ell}$ or $K^{\pm \ell/2}$ defines an orthogonal structure on $V$ and $\Lambda^nV$ is trivial so it has structure group $SO(n,\C)$. 

 The Higgs field must be symmetric with respect to this orthogonal structure. We set:
 \begin{equation}
 \Phi=\begin{pmatrix}
 0 & 1 & 0 &\dots & & 0\\
a_2 & 0 & 1 & \dots & & 0\\
a_3 & a_2 & 0 & 1& \dots & 0\\
\vdots & & &\ddots & &\vdots \\
a_{n-1}& & & & \ddots & 1\\
a_n& a_{n-1} & \dots & a_3 & a_2& 0
\end{pmatrix}
 \label{sym}
 \end{equation}
 where $a_i\in H^0(\Sigma, K^i)$. 
 Then the higher Teichm\"uller space is the  space  $H^0(\Sigma,K^2)\oplus H^0(\Sigma,K^3)\oplus \cdots \oplus H^0(\Sigma, K^n)$ of differentials of degree $2$ to $n$.

The analogy in the case of $SL(\infty,\R)$ was initiated in \cite{Hit2}. We interpret $SU(\infty)$ as the symplectic diffeomorphisms of $S^2$ so a principal $SU(\infty)$-bundle with connection is a 2-sphere bundle $\pi:M\rightarrow\Sigma$ with a symplectic connection. We take this to be the projective bundle $\PP(K^{-1/2}\oplus K^{1/2})\cong \PP(1\oplus K)$. The group $SO(\infty)\subset SU(\infty)$ is taken to be the subgroup which commutes with reflection about an equator in $S^2$, so $M$ has a fibre-preserving involution $\tau$ which is anti-symplectic and fixes a circle bundle. The Higgs field is interpreted as a holomorphic section of the pull back $\pi^*K$ on $M$. Such a section defines a map to the total space of $K$, i.e. the cotangent bundle, and then it is the pullback of the tautological section $s$.The Higgs field should be invariant by $\tau$  which means that $M/\tau$, the disc bundle, maps to $T^*\Sigma$. 

As to the analogue of the equation $F_A+[\Phi,\Phi^*]=0$, this becomes the equation for a hyperk\"ahler structure on $M$, though necessarily with a singularity which is called a {\it fold}.  Recall that a hyperk\"ahler 4-manifold is determined by three closed symplectic 2-forms $\omega_1,\omega_2,\omega_3$ such that the exterior products satisfy
$$\omega_1^2=\omega_2^2=\omega_3^2,\qquad \omega_1\omega_2=\omega_2\omega_3=\omega_3\omega_1 = 0.$$
Then these are K\"ahler forms for complex structures $I,J,K$ satisfying the quaternionic identities. Here we take $\omega_2+i\omega_3$ to be the pullback of the canonical holomorphic symplectic form on $T^*\Sigma$ and $\omega_1$ a closed 2-form which defines the $SU(\infty)$ connection. We have $\tau^*\omega_1=-\omega_1$ and $\tau^*(\omega_2+i\omega_3)=(\omega_2+i\omega_3)$. 
\begin{ex}
The standard model for this is to take the Higgs bundle for the hyperbolic metric, where all terms belong to finite-dimensional subgroups $SO(2), SL(2,\R)$  of the $n=\infty$ groups and then, 
 using the canonical basis $\omega,\theta,\eta$ for the unit circle bundle,  we obtain the hyperk\"ahler metric on the image in $T^*\Sigma$, namely the unit disc bundle of the hyperbolic metric, as
\begin{equation}
\omega_2+i\omega_3 = d(\sqrt{1-u^2}(\omega+i\theta))\qquad \omega_1=d(u\wedge \eta)
\label{hypk} 
\end{equation} 
for $u\in [0,1]$.  At $u=0$ we have a singularity in  $\omega_1$, but this is caused by the folding map $M\rightarrow M/\tau\rightarrow T^*\Sigma$. The form   $\omega_1$ restricts to a smooth symplectic form on the fibre $S^2$. 
\end{ex}

For hyperk\"ahler metrics of this type we have an analogue of the coefficients of the characteristic polynomial, or more precisely the sections of $K^m$ given by $\tr \Phi^m$. We have $s^m$ a holomorphic section of $\pi^*K^m$, and  integrating $\omega_1s^m$ over the 2-sphere fibres of $\pi:M\rightarrow \Sigma$ yields a holomorphic section of $K^m$ on $\Sigma$. This gives zero in the standard model. It is this observation that suggests the existence of such hyperk\"ahler structures for any sum of holomorphic sections, or, in relation to the previous section to a CR function $h$ on $U$. For small values of $h$ this was achieved in Biquard's paper \cite{B1}.  

For a CR function $h$ given by a section locally written $f(z)dz^m$ of $K^m$ the first order variation of the hyperk\"ahler forms is given in \cite{Hit2} as $\dot\omega_1=0$ and 
$$\dot\omega_2+i\dot\omega_3=d(f(z)dz^m(\bar wd\bar z)^{m-1}y^{2m-2}(dz d\bar z)^{-(m-1)})$$
where more invariantly $dzd\bar z/y^2$ is the hyperbolic metric and $wdz$ the canonical 1-form. In the terminology of equation (\ref{hypk}) we have $wdz=\sqrt{1-u^2}(\omega+i\theta)$ and $h=s^{-m}f(z)dz^m$ giving 
\begin{equation} \label{defhk} 
\dot\omega_2+i\dot\omega_3=d(h(\sqrt{1-u^2})^{m-1}(\omega+i\theta)).
\end{equation}

Each solution, at least for small deformations,  defines a deformed disc bundle in $T^*\Sigma$ which is convex and contains the origin and thus defines a Finsler metric of negative curvature, whose geodesics we conjecture should play a significant role.

\begin{rmk}
It is important to note  (Theorem 3 in \cite{Hit2}) that the only Riemannian metrics in this family are hyperbolic. Nevertheless, the flat symplectic connection for a Riemannian surface could still be equivalent to that of a different Finsler surface -- for example Bonahon \cite{Bon} has shown that the same marked length spectrum can be realized by a Riemannian and a  non-Riemannian metric structure.  
\end{rmk} 

\subsection{The hypersymplectic extension} 
Pursuing the analogy further comes up against the serious problem of the complex group -- there is no obvious complexification of $SL(\infty,\R)$ that fits into the picture.   Donaldson's substitutes for the complexification of a group of symplectic diffeomorphisms may  be relevant and certainly some aspects of \cite{Don1},\cite{Don2} seem close to the issue of folded hyperk\"ahler metrics, but here we follow a different line of thought. 
    
   Take the hyperk\"ahler example  (\ref{hypk}) and put $u=i\sinh (t)$ then  we obtain  closed forms $\varphi_1,\varphi_2,\varphi_3$ given by $\omega_1=i\varphi_1$, $\omega_2=\varphi_2$, $\omega_3= \varphi_3$ satisfying 
   \begin{equation}
   -\varphi_1^2=  \varphi_2^2=\varphi_3^2,\quad \varphi_1\varphi_2=\varphi_2\varphi_3=\varphi_3\varphi_1 = 0 
   \label{hyps}
   \end{equation}
and moreover $$\varphi_2-\varphi_1=d(\cosh (t) \omega-\sinh (t)\eta)$$ which we recognize as the symplectic connection for a hyperbolic metric.  

A 4-manifold with closed symplectic 2-forms satisfying (\ref{hyps}) is called {\it hypersymplectic} \cite{HitX},\cite{DS}. Whereas a hyperk\"ahler 4-manifold has a Riemannian metric with holonomy $SU(2)$ a hypersymplectic 4-manifold has a signature $(2,2)$ metric with holonomy $SU(1,1)$. 
The characteristic property is that 
$(\varphi_1+\cos\theta \varphi_2+\sin \theta \varphi_3)^2=0$ and since each such form is closed and locally decomposable, for fixed $\theta$ its annihilator defines a foliation by isotropic surfaces. Putting $\theta=-\pi/2$ in the example gives the foliation for the flat symplectic connection. This suggests the relationship holds in more generality. 

Biquard's hyperk\"ahler manifolds are real analytic and so each one has an analytic  continuation to an external neighbourhood of the boundary of the disc bundle. In fact more is true. Given a real analytic 3-manifold $U$ together with closed 2-forms $\beta_2,\beta_3$ whose kernel generates a contact distribution, then on $(-\epsilon,\epsilon)\times U$ there exists a real analytic folded hyperk\"ahler metric with $\omega_2,\omega_3$ restricting to $\beta_2,\beta_3$ and with involution $\tau$ satisfying $\tau^*\omega_1=-\omega_1$ and leaving $\omega_2,\omega_3$ fixed. Thus Biquard's globally defined metrics are uniquely determined by the boundary data on the fold. This is Theorem 3 in \cite{B1} and Theorem 1 in \cite{Hit2}. However Remark 2 following the latter proof observes that there is also a natural local {\it hypersymplectic} extension with the same property.  Then $(-\epsilon,\epsilon)\times U$
is a hypersymplectic manifold with a $\tau$-invariant map to $T^*\Sigma$ whose image is the analytic continuation of Biquard's metric.   It is a $(-\epsilon,\epsilon)\times S^1$ bundle over $\Sigma$.

Since $\tau^*\varphi_1=-\varphi_1$, $\varphi_1$  vanishes when restricted as a form to  the fixed point set $U$, so $\alpha=\varphi_2- \varphi_1$ is a closed decomposable 2-form which restricts on $U$ to $\varphi_2$. But $\varphi_2$ is  the analytic continuation of the canonical symplectic form on $T^*\Sigma$ which exists on the whole cotangent bundle and so is still the canonical form. Thus a leaf of the foliation  defined by $\varphi_2- \varphi_1$ intersects $U$ in an orbit of the geodesic flow for the Finsler metric defined by $U$. The symplectic connection  can therefore locally be described by parallel translation of geodesics but until we can step beyond the local situation  we can't make any statement about the global holonomy in general.

Nevertheless, putting $u=i\sinh(t)$ in the variation \ref{defhk} for the analytic continuation, we obtain 
$$\dot\varphi_2+i\dot\varphi_3=d(h\cosh^{m-1}(t)(\omega+i\theta))$$
and since $\dot\varphi_1=0$ this means the deformation of the symplectic connection $\dot\varphi_1-\dot\varphi_2$ is the CR deformation of Proposition \ref{CRprop}. 

\subsection{Involutions} 
A further piece of evidence for the role of the hypersymplectic structure comes when we look at Corlette's approach to Higgs bundles, beginning with the flat connection. For a real group $G$ with maximal compact subgroup $H$ and an irreducible  representation of $\pi_1(\Sigma)$ in $G$,  the theorem in \cite{Cor} produces an equivariant harmonic map $f:\tilde \Sigma \rightarrow G/H$. The subgroup $H$ is 
the fixed point set of an involution $i$.

Consider the Hamiltonian diffeomorphisms  of $S^1\times \R$ with symplectic form $\varphi=dt\wedge d\theta$ and 
define the antisymplectic involution  $\sigma:S^1\times \R\rightarrow S^1\times \R$  by $\sigma(\theta, t)=(\theta,-t)$. Putting $G=\Ham(S^1\times \R)$ we have the automorphism $i$ of $G$ given by $i(g)=\sigma g \sigma^{-1}$, with stabilizer $H$. Formally, this defines a symmetric space $G/H$.

\begin{rmks} 

\noindent 1. If we write $G=SL(\infty,\R)$ it would be natural to set $H=SO(\infty)$ but this is not the same group as our other usage which is the symplectic diffeomorphisms of $S^2$ which  commute with reflection in an equator -- one of the other limitations of drawing an analogy with the finite-dimensional case.

\noindent 2. Now $PSL(2,\R)\subset G$ and $SO(2)\subset H$ so each point $x$ of the hyperbolic plane $PSL(2,\R)/SO(2)\subset G/H$ defines an involution on the space of oriented geodesics -- the fixed point set consisting of the geodesics through $x$.  In classical geometric terms we can use the Beltrami model of the unit disc $D$ as the interior of the conic $x_1^2+x_2^2-x_3^2=0$ in $\RP^2$ with geodesics segments of projective lines. If $(x,x)<0$ represents a point in $D$, an exterior point $y$ (where $(y,y)>0$) has a polar line which intersects $D$ in a geodesic. Then the involution on (unoriented) geodesics is just the reflection using the indefinite inner product
$$\sigma_x(y)=y-2\frac{(x,y)}{(x,x)}y.$$
\end{rmks}

Suppose now we have a representation of $\pi_1(\Sigma)$ in $G=SL(\infty,\R)$ and an equivariant map $f:\tilde \Sigma\rightarrow G/H$. Then for $x\in \tilde \Sigma$, $f(x)$ defines an involution $\sigma_x$.

 Take the product $M=\tilde \Sigma\times S^1\times \R$ so that $\varphi=dt\wedge d\theta$ defines the trivial symplectic connection and define 
$\tau(x,y)=(x,\sigma_x(y)).$
Set $2\varphi_1=\tau^*\varphi-\varphi, 2\varphi_2=\tau^*\varphi+\varphi$ and assume these forms are symplectic. By definition they satisfy the relations 
$$\varphi_1\varphi_2=0,\qquad \varphi_1^2=-\varphi_2^2$$
and $\varphi_2-\varphi_1=\varphi$ is the flat symplectic connection corresponding to the product.

We express  the tangent bundle as $TM\cong H\oplus V$, the horizontal/vertical decomposition. Since $\tau$ is anti-symplectic on the fibres $\varphi_1=(\varphi+\tau^*\varphi)/2$ lies in the subbundle $\Lambda^2H^*\oplus (H^*\otimes V^*)\subset \Lambda^2T^*M$. Denote by $\rho$ the component in 
$H^*\otimes V^*$. A conformal structure on $\Sigma$ defines a Hodge star operator on 1-forms so acting on the $H^*\cong \pi^*T^*\Sigma$ factor we have a 2-form $\ast\rho$. Since $\alpha\wedge\ast\beta$ is symmetric we have $\rho\wedge \ast\rho=0$. And since $H^*$ is 2-dimensional it follows that 
$\varphi_1\wedge \ast\rho=0$.

Now $\varphi_1^2=\rho^2$ since $\Lambda^pH^*=0$ for $p>2$ and $\rho^2=(\ast\rho)^2$ so we have the algebraic relations
$$\varphi_1^2=(\ast \rho)^2,\qquad \varphi_1\wedge \ast\rho=0.$$
Furthermore $\varphi_2=\varphi_1+2\varphi$ so that $\varphi_2\wedge \ast\rho=0$ since $\varphi$ lies in the $\Lambda^2V^*$ component. If we write $\varphi_3=\ast\rho$ then we have 
$$\varphi_3^2=\varphi_1^2=-\varphi_2^2\qquad \varphi_1\varphi_2=\varphi_2\varphi_3=\varphi_3\varphi_1=0.$$

Now $d\varphi_1=0=d\varphi_2$. If $d\varphi_3=0$ these are the defining equations for a hypersymplectic structure on $M$, assuming the nondegeneracy of the forms.

A symmetric space  $G/H$ is a space of involutions, all of which are conjugate to a fixed one, which gives an embedding in the group $G$. Then $f:\tilde\Sigma \rightarrow G/H$ can be defined by a map $g:\tilde \Sigma\rightarrow G$. In finite dimensions this map is harmonic if 
$d(\ast g^{-1}dg)=0.$ We adopt this definition in our infinite-dimensional situation, and then we have

\begin{prp} \label{harm}The form $\varphi_3=\ast \rho$ is closed if and only if the map $f:\tilde \Sigma\rightarrow G/H$ is harmonic.
\end{prp}
\begin{prf} 
The map $g$  on the universal covering is the gauge transformation which relates the two trivializations -- one using the flat connection  $\varphi$, and the other using  $\tau^*\varphi$, and so $g^{-1}dg$ is the difference of the two connections. But the difference of two symplectic connections  is a section of $H^*\otimes V$, and using the isomorphism $V^*\cong V$ as in Proposition \ref{inf} this  section is the component of $\varphi+\tau^*\varphi$ in $H^*\otimes V^*$ or $2\rho$. 

Thus if we take the definition of a harmonic map to be $d(\ast g^{-1}dg)=0$, this is equivalent to $d(\ast\rho)=0$. 
\end{prf}
The proposition provides further motivation for regarding the hypersymplectic extension as playing a key role.
\subsection{Quadratic differentials} 
In this final section we evaluate the hypersymplectic structure based on a harmonic map  $f: \tilde \Sigma\rightarrow PSL(2,\R)/SO(2)\subset G/H$. This is the basis of the Higgs bundle description of classical Teichm\"uller space by quadratic differentials $q$, and is the first family of CR deformations. In this case the image of $\pi_1(\Sigma)$ is an isomorphic subgroup $\Gamma\subset SL(2,\R)$ and the equivariant map descends to  a  diffeomorphism $\Sigma\rightarrow \Gamma\backslash SL(2,\R)/SO(2)$.   The original existence  is due to
 \cite{EE}, though here we follow   Section 11 of \cite{Hit0}. 
 
 In the Higgs bundle interpretation, the flat $SL(2,\R)$ bundle is given by  $V=K^{-1/2}\oplus K^{1/2}$ with $\Phi(u,v)=(v,qu)$
 and a hermitian structure on $K^{1/2}$, or equivalently a hermitian form $(\,\,,\,\,)$ on $K$. Then we have 
 \begin{prp} The equivariant harmonic map $f$ defines a folded hypersymplectic structure on the 4-manifold $M=\{(u, x)\in K\times \R: (u,u)=x^2\}$. In  a local coordinate $z$ on $\Sigma$  it is defined by the three closed two-forms $\varphi_1,\varphi_2,\varphi_3$:
 $$ \varphi_1=-2d(\sinh (t)(d\phi-i\alpha))$$
 $$\varphi_2+i\varphi_3=d(\cosh (t) (h^{-1/2}ae^{i\phi}-e^{-i\phi}h^{1/2}) dz)$$
 where $x=\sinh t, q=a(z)dz^2$, the Hermitian form on the tangent bundle is $hdzd\bar z$ and $\alpha$ is the connection form $(h^{-1}\partial h -h^{-1}\bar\partial h)/2$. 
  \end{prp} 
 \begin{prf} 
 
 %An equivariant harmonic map in this case is a map from a surface $\Sigma$ with a conformal structure to a hyperbolic surface $S$ of the same %genus. The harmonicity condition implies that the pull-back of the 
%metric has a holomorphic $(2,0)$ component which is a quadratic differential $q\in H^0(\Sigma, K^2)$. Conversely, given $q$, there exists a hermitian metric $h\in \Omega^{1,1}(\Sigma)$ such that 
%\begin{equation}
%\hat h=q+(h+\frac{q\bar q}{h})+\bar q
%\label{metric}
%\end{equation} 
%is a metric of constant curvature $-4$. 

% and we have a $U(1)$ connection $A$  and equations $F_A+[\Phi,\Phi^*]=0$ giving the flat connection $\nabla_A+\Phi+\Phi^*$. The Hermitian form on $K^{1/2}$ is equivalent to a Hermitian metric $hdzd\bar z$ in terms of a local holomorphic coordinate $z$ on $\Sigma$. 
 The rank $2$ Higgs bundle is $V=K^{-1/2}\oplus K^{1/2}$ with $\Phi(u,v)=(v,qu)$ and in the given notation $(h^{-1/4}dz^{-1/2}, h^{1/4}dz^{1/2})$ is a local unitary frame for $V$. 
If $(u,v)$ are the coefficients of a local section  in this frame, $V$ has a real structure $(u,v)\mapsto (\bar v, \bar u)$ which is preserved by the flat connection $\nabla_A+\Phi+\Phi^*$, giving holonomy in $SL(2,\R)$. 

This is the action of $SL(2,\R)$ on $\R^2$ and we need the action on the symplectic surface $S^1\times \R$, which we identify as a coadjoint orbit. This means we take the Lie algebra bundle $\End_0V=K^{-1}\oplus 1\oplus K=S^2V$. This has an induced real structure and a Killing form. Using coefficients $(u_{-1},u_0,u_1)$ with respect to the frame $(h^{-1/2}dz^{-1},1,h^{1/2}dz)$ of $K^{-1}\oplus 1\oplus K$ the coadjoint orbit is $u_{-1}u_1-u_0^2=1$ and the real structure $(u_{-1},u_0,u_1)\mapsto(\bar u_1,\bar u_0,\bar u_{-1})$. 

We  parametrize the real surface with $u_1=\cosh (t)e^{i\phi}, u_0=\sinh (t)$, then this $S^1\times \R$ bundle over $\Sigma$ is invariantly the description  $M=\{(u, x)\in K\times \R: (u,u)=x^2\}$ in the statement of the proposition. It is a symplectic bundle and in local coordinates the form is  $du_1\wedge d\bar u_1/2iu_0=-\cosh (t) dt\wedge d\phi$. 
The  involution $\tau$ on $M$ determined by the reduction to $SO(2)$ is $t\mapsto -t$.
 
 Spelling out the connection $\nabla_A +\Phi+\Phi^*$ in coordinates, we see that   a local covariant constant section of $M\rightarrow \Sigma$ is given by the vanishing of   the complex form $$\beta=du_1+\alpha u_1+(h^{1/2}dz+\bar ah^{-1/2}d\bar z)u_0=\beta_1+u_0\beta_2$$
 where $q=adz^2$ and $\alpha=h^{-1}(\partial h/\partial z-\partial h/\partial \bar z)/2$ is the connection form for $\nabla_A$ in the unitary basis of $K^*$, the tangent bundle.  Then $\varphi=\beta\wedge \bar\beta/2iu_0$ is a closed real $2$-form restricting to the symplectic form on the fibre. 
 
  We have 
 $$u_0^{-1}\beta\wedge\bar\beta=u_0^{-1}\beta_1\wedge\bar\beta_1+u_0\beta_2\wedge\bar\beta_2+\beta_1\wedge\bar\beta_2+\beta_2\wedge\bar\beta_1.$$
  In the $t,\phi$ coordinates  we see that the odd part under $u_0\mapsto -u_0$, namely  
 $\varphi_1=\frac{1}{2}(\tau^*\varphi-\varphi)$ is $$2\cosh (t) (d\phi-i\alpha)\wedge dt-\sinh (t) (h-h^{-1}\vert a\vert^2)idz\wedge d\bar z 
 = -2d(\sinh t(d\phi-i\alpha))$$
using $d\alpha =  F=h-h^{-1}\vert a\vert^2$  from the equation $F_A+[\Phi,\Phi^*]=0$. The 2-form $\varphi_1$  is nondegenerate if $t>0$ and $h-h^{-1}\vert a\vert^2\ne 0$. The latter holds by the maximum principle as in \cite{Hit0}.

For the even part note that $d\bar\beta_2=d(h^{1/2}d\bar z+ ah^{-1/2}dz)=\alpha\wedge(h^{1/2}d\bar z+ ah^{-1/2}dz)$ since $a$ is holomorphic, so 
$\beta_1\wedge \bar\beta_2=d(u_1\bar\beta_2)$ and 
$$\varphi_2=\frac{1}{2}(\tau^*\varphi+\varphi)=\frac{1}{2i}d(u_1\bar\beta_2-\bar u_1\beta_2).$$
Now in $\beta_1\wedge \bar\beta_2$ the second factor lies in $H^*$ and $\ast \bar\beta_2=ih^{1/2}d\bar z-iah^{-1/2}dz$ so, as above, we obtain  
$$\varphi_3=\frac{1}{2i}d(u_1\ast\bar\beta_2-\bar u_1\ast\beta_2)$$
(this is clearly closed, as implied by Proposition \ref{harm}). We can now write  
$$\varphi_2+i\varphi_3=d(\cosh (t) (h^{-1/2}ae^{i\phi}-e^{-i\phi}h^{1/2}) dz)$$
\end{prf} 
%is a section of $V^*\otimes H^*$ and 
%$$\ast \beta_2=\ast (h^{1/2}dz+\bar ah^{-1/2}d\bar z)=i(-h^{1/2}dz+\bar ah^{-1/2}d\bar z)$$
%so that $\beta_2+i\ast \beta_2=  2h^{1/2}dz$ and $\bar\beta_2+i\ast\bar\beta_2= 2ah^{1/2}dz.$
 Our motivation for introducing the hypersymplectic manifold was to generate a symplectic connection for the Finsler structure defined by the fold. In this case the fold is the circle bundle over $\Sigma$ given by $t=0$ -- the  unit circle for the hermitian form $h$. Then we have 
 \begin{prp} The symplectic connection defined by the above hypersymplectic structure is equivalent to the connection defined in Theorem \ref{main} for the hyperbolic surface $\Gamma\backslash SL(2,\R)/SO(2)$.
\end{prp} 
\begin{prf} 
Consider the real 2-form $\varphi_2=d(u_1\bar\beta_2-\bar u_1\beta_2)/2i$ restricted to $t=0$,  which is $\vert u_1\vert=1$. Then $u_1\bar\beta_2-\bar u_1\beta_2$ is a section of $\pi^*T^*\Sigma$ and maps $U$ to the total space of the cotangent bundle. With $u_1=e^{i\phi}$, the image preserves the metric $(\beta_2\otimes \bar\beta_2+\bar\beta_2\otimes \beta_2)/2$. Since $\beta_2= h^{1/2}d z+ \bar ah^{-1/2}d\bar z$ the metric  is 
$$adz^2+\left(h+\frac{a\bar a}{h}\right)dzd\bar z+\bar a d\bar z^2$$
and, as in  \cite{Hit0}, this is the hyperbolic metric corresponding to $\Gamma\subset SL(2,\R)$ with $\beta_2$ a local unit $(1,0)$ form. Thus $U$ maps to the unit circle bundle $\hat U$ of this metric.

The map $U\rightarrow \hat U\subset T^*\Sigma$ means that we can identify $(e^{i\phi}\bar\beta_2-e^{-i\phi}\beta_2)/2i$ with the restriction of the canonical one-form $\omega$ on $T^*\Sigma$. From $d\bar\beta_2=\alpha\wedge\bar\beta_2$ we obtain
$$d\omega = d(e^{i\phi}\bar\beta_2-e^{-i\phi}\beta_2)/2i=(d\phi-i\alpha)\wedge(e^{i\phi}\bar\beta_2+e^{-i\phi}\beta_2)/2.$$
  The Hodge star $\hat\ast$ for this metric takes $\bar\beta_2$ to $i\bar\beta_2$ and so $\hat\ast \,\omega=(e^{i\phi}\bar\beta_2+e^{-i\phi}\beta_2)/2 =\theta$ in the notation of Theorem 2. Setting $\eta=d\phi-i\alpha$, the connection form gives the isomorphism with $2d(\cosh(t) \omega-\sinh(t) \eta)$, the connection for constant curvature. 
\end{prf}

\end{document}